\documentclass{amsart}
\usepackage{amsfonts}
\usepackage{amsmath}
\usepackage{amssymb}
\usepackage{amscd}
\usepackage{amstext}
\usepackage{amsthm}
\usepackage{enumerate}
\usepackage{color}

\theoremstyle{plain}
\newtheorem{lem}{Lemma}[subsection]
\newtheorem{prop}[lem]{Proposition}
\newtheorem{thm}[lem]{Theorem}

\newtheorem{fact}[lem]{Fact}

\newtheorem{mainthm}{Theorem}

\theoremstyle{definition}
\newtheorem{defn}[lem]{Definition}

\theoremstyle{remark}
\newtheorem{ex}[lem]{Example}
\newtheorem{rem}[lem]{Remark}

\begin{document}
\title{Reduction of generalized Calabi-Yau structures}
\author{Yasufumi Nitta}
\date{}
\maketitle
\thispagestyle{empty}

\begin{abstract}
A generalized Calabi-Yau structure is a geometrical structure on a manifold which 
generalizes both the concept of the Calabi-Yau structure and that of the symplectic one. 
In view of a result of Lin and Tolman in generalized complex cases, we introduce in this paper 
the notion of a generalized moment map for a compact Lie group action 
on a generalized Calabi-Yau manifold 
and construct a reduced generalized Calabi-Yau structure on 
the reduced space. As an application, we show some relationship between generalized moment maps and 
the Bergman kernels, and prove the Duistermaat-Heckman formula 
for a torus action on a generalized Calabi-Yau manifold.
\end{abstract}

\section{Introduction}
Generalized Calabi-Yau structures introduced by Hitchin \cite{Hi} were developed by Gualtieri \cite{Gua} as a special case of 
generalized complex structures. It is a geometrical structure defined by a differential form, which generalizes both the concept of the Calabi-Yau 
structure - a non vanishing holomorphic form of the top degree - and that of the symplectic structure. In this paper, we consider 
a compact Lie group action on a generalized Calabi-Yau manifold. 

A compact Lie group action on a generalized complex manifold was studied by Lin and Tolman in \cite{Lin1} and \cite{Lin2}. In \cite{Lin1}, 
they introduced a notion of generalized moment maps for a compact Lie group action on a generalized complex manifold by generalizing the notion of 
moment maps for a compact Lie group action on a symplectic manifold. Using this definition, they constructed a generalized complex structure 
on the reduced space, which is natural up to a transformation by an exact $B$-field. 

In the present paper, we apply the definition of a generalized moment map to a compact Lie group action on a generalized Calabi-Yau manifold, 
and construct a generalized Calabi-Yau structure on the reduced space. Moreover, we shall show that the reduced generalized Calabi-Yau structure 
is unique and has the same type as the original generalized Calabi-Yau structure (cf. Section 3).

\begin{mainthm}\label{main1}
Let a compact Lie group G act on a generalized Calabi-Yau manifold $(M, \varphi)$ in a Hamiltonian way with a generalized moment map $\mu: M \longrightarrow \frak{g}^*.$ 
If G acts freely on $\mu^{-1}(0)$, then the quotient space $M_0 = \mu^{-1}(0)/G$ is a smooth manifold, and inherits a unique generalized Calabi-Yau structure $\widetilde \varphi$ which satisfies
\begin{eqnarray*}
p_{0}^{*} \widetilde \varphi = i_{0}^{*} \varphi, 
\end{eqnarray*}
where $i_{0} : \mu^{-1}(0) \longrightarrow M$ is the inclusion and $p_{0}: \mu^{-1}(0) \longrightarrow M_0$ is the natural projection. 
Moreover, for each $p \in \mu^{-1}(0)$, 
\begin{eqnarray*}
{\rm type}(\varphi_p) = {\rm type}(\widetilde \varphi_{[p]}). 
\end{eqnarray*}
\end{mainthm}

The detailed definitions of the theorem are in Section 3. In particular, in the case that the generalized Calabi-Yau structure 
is induced by a symplectic structure, the reduced form is induced by the reduced symplectic form. In addition we construct an example of a 
Hamiltonian action on a generalized Calabi-Yau structure which is not induced by either a symplectic structure or a Calabi-Yau one. 
We then show some relationship between generalized moment maps and Bergman kernels (cf. $\bold{Example}$ \ref{ex1} and \ref{ex2} in Section 3). 

We next consider that a generalized Calabi-Yau structure $\varphi$ on a connected manifold $M$ which has constant type $k$. 
Then there exists a natural volume form $dm = \frac{(\sqrt{-1})^n}{2^{n-k}} \langle \varphi , \bar \varphi \rangle$ 
defined by $\varphi$, which generalizes the Liouville form on a symplectic manifold. Indeed, if $\varphi$ is 
a generalized Calabi-Yau structure induced by a symplectic structure $\omega$, 
then $dm$ coincides with the Liouville form for the symplectic structure $\omega$.
Further by assuming that a compact torus $T$ acts on $M$ effectively.
Under the assumptions, we shall show the Duistermaat-Heckman formula for the volume form $dm$ (cf. Section 4). 

\begin{mainthm}\label{main2}
Let $(M, \varphi)$ be a $2n$-dimensional connected generalized Calabi-Yau manifold which has constant type $k$, and suppose that compact $l$-torus $T$ 
acts on $M$ effectively and in a Hamiltonian way. In addition, we assume that the generalized moment map $\mu$ is proper. 
Then the pushforward $\mu_{*}(dm)$ of the natural volume form $dm$ under $\mu$ is absolutely continuous with respect to the 
Lebesgue measure on $\frak{t}^{*}$ and the Radon-Nikodym derivative $f$ can be written by 
\[
f(a) = \int_{M_{a}}dm_{a} = {\rm vol}(M_{a}) 
\]
for each regular value $a \in \frak{t}^{*}$ of $\mu$, and 
$dm_{a}$ denotes the measure defined by the natural volume form on the reduced space $M_{a} = \mu^{-1}(a)/T$. 
\end{mainthm}

This paper is organized as follows. In Section 2 we introduce background materials and the definition of generalized Calabi-Yau structures. 
In Section 3 we define the notion of generalized moment maps for a Lie group action on a generalized Calabi-Yau manifold, 
and construct a generalized Calabi-Yau structure on the reduced space. In addition, we discuss some relations 
between generalized moment maps and Bergman kernels. At last Section, we proved  the Duistermaat-Heckman formula 
for a Hamiltonian torus action on a generalized Calabi-Yau manifold. 
%%%%%%%%%%%%%%%%%%%%%%%%%%%%%%%%%%%%%%%%%%%%%%%%%%%%%%%%%%%%%%%%%%%%%%%%%%%%%%%%%%%%%%%%%%%%%%%%%%%%%%%%%%%%%%%%%%%%%%%%%%%%%
%%%%%%%%%%%%%%%%%%%%%%%%%%%%%%%%%%%%%%%%%%%%%%%%%%%%%%%%%%%%%%%%%%%%%%%%%%%%%%%%%%%%%%%%%%%%%%%%%%%%%%%%%%%%%%%%%%%%%%%%%%%%%
\section{Generalized Calabi-Yau structures}
In this section we recall the definition of generalized Calabi-Yau structures. 
For the detail, see \cite{Gua} and \cite{Hi}.

\subsection{Clifford algebras and the spin representation}
Let $V$ be a real vector space of dimension $n$, and $V^{*}$ be the dual space of $V$. Then the direct sum $V \oplus V^{*}$ 
admits a natural indefinite metric of signature ($n$,\ $n$) defined by 
\begin{eqnarray*}
( X+\alpha,\ Y + \beta ) = \frac{1}{2}\left( \beta (X) + \alpha (Y) \right) 
\end{eqnarray*}
for $X+\alpha,\ Y+\beta \in V \oplus V^*$. 
Let ${\rm T}(V \oplus V^*) = \oplus_{p=0}^{\infty} \left(\otimes^{p}(V \oplus V^{*}) \right)$ be the tensor algebra of $V \oplus V^{*}$, and 
define $\mathcal{I}$ to be the two-sided ideal generated by 
$\{(X+\alpha)\otimes(X+\alpha) - ( X+\alpha,\ X+\alpha ) \big|\ X+\alpha \in V \oplus V^* \}\}$. 
Then we call the quotient algebra
\begin{eqnarray*}
{\rm CL}(V \oplus V^{*}) = {\rm T}(V \oplus V^{*})/\mathcal{I}
\end{eqnarray*}
the Clifford algebra of $V \oplus V^{*}$. 
For each $E, F \in {\rm CL}(V \oplus V^*)$, $E\cdot F$ denotes the multiplication induced by the tensor product. 

Consider the exterior algebra $\wedge^{*}V^{*}$ and a linear mapping $V \oplus V^{*} \longrightarrow {\rm End}(\wedge^{*}V^{*})$ 
defined by 
\begin{eqnarray*}
(X+\alpha) \cdot \varphi = \iota_{X}\varphi + \alpha \wedge \varphi. 
\end{eqnarray*}
Then we have 
\begin{eqnarray*}
(X + \alpha)^2 \cdot \varphi &=& \iota_{X}(\alpha \wedge \varphi) + \alpha \wedge \iota_{X}\varphi \\
                             &=& (\iota_{X}\alpha)\varphi \\
                             &=& ( X + \alpha,\ X + \alpha ) \varphi,
\end{eqnarray*}
so it can be extended to a representation of the Clifford algebra ${\rm CL}(V \oplus V^{*}) \longrightarrow {\rm End}(\wedge^{*}V^{*})$. 
This is called the spin representation, and a element $\varphi \in \wedge^{*}V^{*}$ is called a spinor. 

We define ${\rm Pin}(V \oplus V^*)$ and ${\rm Spin}(V \oplus V^*)$, subgroups of the group consists of invertible elements of 
${\rm CL}(V \oplus V^{*})$ by 
\begin{eqnarray*}
{\rm Pin}(V \oplus V^*) = \{E_1 \cdots E_k \big|\ k \in \mathbb{N}\cup \{0\},\ ( E_i , E_i ) = \pm 1\}, \\
{\rm Spin}(V \oplus V^*) = \{E_1 \cdots E_{2k} \big|\ k \in \mathbb{N}\cup \{0\},\ ( E_i , E_i ) = \pm 1\}.
\end{eqnarray*}
we call ${\rm Pin}(V \oplus V^*)$ the pin group, and ${\rm Spin}(V \oplus V^*)$ the spin group. The following proposition says a 
geometrical meaning of the pin and spin group. 

\begin{prop}[\cite{Che}, \cite{Gua}]
The pin group and the spin group have following short exact sequences. 
\begin{eqnarray*}
1 \longrightarrow \mathbb{Z}/2\mathbb{Z} \longrightarrow {\rm Pin}(V \oplus V^*) \longrightarrow {\rm O}(V \oplus V^{*}) \longrightarrow 1 \\
1 \longrightarrow \mathbb{Z}/2\mathbb{Z} \longrightarrow {\rm Spin}(V \oplus V^*) \longrightarrow {\rm SO}(V \oplus V^{*}) \longrightarrow 1 \\
\end{eqnarray*}
\end{prop}

Let ${\rm Spin}_{0}(V \oplus V^*)$ denote the identity component of ${\rm Spin}(V \oplus V^*)$. 
Then $\wedge^{*}V^{*}$ has a ${\rm Spin}_{0}(V \oplus V^*)$-invariant bilinear form defined by 
\begin{eqnarray*}
\langle \varphi, \psi \rangle = (\sigma(\varphi) \wedge \psi )_{n}, 
\end{eqnarray*}
where $()_{n}$ indicates taking the $n$-th degree component of the form, and $\sigma : \wedge^{*}V^* \longrightarrow \wedge^{*}V^*$ is an anti-homomorphism on $\wedge^{*}V^{*}$
defined by
\begin{eqnarray*}
\sigma(\varphi_1 \wedge \cdots \wedge \varphi_k) = \varphi_k \wedge \cdots \wedge \varphi_1
\end{eqnarray*} 
for each $\varphi_1, \cdots, \varphi_k \in \wedge^{1}V^{*}$. 

\subsection{Pure spinors and generalized Calabi-Yau structures on a vector space}
Given a spinor $\varphi \in \wedge^{*}V^{*}$, we define the annihilator of $\varphi$ by 
\begin{eqnarray*}
E_{\varphi} = \{ X + \alpha \in V \oplus V^* \big|\ (X + \alpha)\cdot \varphi = 0 \}. 
\end{eqnarray*}
We also define the annihilator $E_{\varphi} \subset (V \oplus V^{*}) \otimes \mathbb{C}$ of a complex spinor $\varphi \in \wedge^{*}V^{*} \otimes \mathbb{C}$ 
in a similar way. Since an element $X + \alpha \in E_{\varphi}$ satisfies 
\begin{eqnarray*}
(X + \alpha, X + \alpha )\varphi &=& (X + \alpha )^{2} \cdot \varphi \\
                                 &=& 0, 
\end{eqnarray*} 
we see that if $\varphi$ is a non-zero spinor or a complex spinor, then $E_{\varphi}$ is isotropic with respect to the 
natural metric on $(V \oplus V^{*}) \otimes \mathbb{C}$. 
In particular, we have $\dim E_{\varphi} \leq n$. 

\begin{defn}
A spinor $\varphi \in \wedge^{*}V^*$ is called pure if $E_{\varphi}$ is maximally isotropic, which means that has the dimension equal to $n$. 
A complex spinor $\varphi \in \wedge^{*}V^* \otimes \mathbb{C}$ with the maximal isotropic subspace $E_{\varphi}$ is called a 
complex pure spinor. 
\end{defn}

\begin{rem}
It is known that if $\varphi \in \wedge^{*}V^*$ is a pure spinor, then $\varphi \in \wedge^{ev/od}V^*$, where 
\begin{eqnarray*}
\wedge^{ev}V^* = \wedge^{0}V^* \oplus \wedge^{2}V^* \oplus \cdots, \\
\wedge^{od}V^* = \wedge^{1}V^* \oplus \wedge^{3}V^* \oplus \cdots, 
\end{eqnarray*}
and $\varphi \in \wedge^{ev/od}V^*$ means that $\varphi$ belongs in either $\wedge^{ev}V^*$ or $\wedge^{od}V^*$. 
\end{rem}

\begin{ex}
The spinor $1 \in \wedge^{0}V^{*}$ is pure, since $E_{1} = V$. 
\end{ex}

\begin{ex}
A non-zero vector $\varphi \in \wedge^{n}V^{*}$ is also pure. The annihilator is $E_{\varphi} = V^{*}$. 
\end{ex}

\begin{ex}
If $\varphi$ is a pure spinor on $V$ and $B$ is a $2$-form, then 
\begin{eqnarray*}
\exp (B)\varphi = \left( 1 + B + \frac{1}{2!}B^{2} + \cdots \right) \wedge \varphi
\end{eqnarray*}
is also pure. The annihilator is $E_{\exp(B)\varphi} = \{ X + \alpha + \iota_{X}B \big|\ X + \alpha \in E_{\varphi} \}$, where $E_{\varphi}$ is the 
annihilator of $\varphi$. 
\end{ex}

Gualtieri shows in his thesis \cite{Gua} that every pure spinor can be written by a 
complex $2$-form and a decomposable complex form as follows. 

\begin{fact}[\cite{Gua}]\label{fact:Gua}
Let $\varphi$ be a complex pure spinor on $V$. 
Then there exists a complex 2-form $B + \sqrt{-1}\omega \in \wedge^{2}V^{*} \otimes \mathbb{C}$ and a complex $k$-form $\Omega$ such that 
\begin{eqnarray*}
\varphi &=& \exp(B + \sqrt{-1}\omega )\Omega. 
\end{eqnarray*}
Moreover, $\Omega$ can be written 
\[
\Omega = \theta^{1}\wedge \cdots \wedge \theta^{k} 
\]
by some 1-forms $\theta^{1}, \cdots,  \theta^{k} \in \wedge^{1}V^{*}\otimes \mathbb{C}$, 
and $\omega$ is nondegenerate on a subspace 
$W = \{X \in V \big|\ \iota_{X}\Omega = 0 \}$. 
\end{fact}
The degree of the form $\Omega$ is called the type of the complex pure spinor $\varphi$ and written by ${\rm type}(\varphi)$. 

Now we give the definition of a generalized Calabi-Yau structure on a real vector space $V$. 

\begin{defn}
Let $V$ be a real vector space of dimension $n = 2m$. 
A generalized Calabi-Yau structure on $V$ is a complex pure spinor $\varphi \in \wedge^{ev/od}V^* \otimes \mathbb{C}$ which 
satisfies that $\langle \varphi, \bar \varphi \rangle \not = 0$.
\end{defn}
Fact \ref{fact:Gua} tells us if there exists a complex pure spinor on $V$ which satisfies 
$\langle \varphi, \bar \varphi \rangle \not = 0$, then $V$ must be even dimensional. 
The condition $\langle \varphi, \bar \varphi \rangle \not = 0$ has the following geometrical meaning.  
\begin{fact}[\cite{Che}]
Let $\varphi$ and $\psi$ be pure spinors. Then they satisfy
$\langle \varphi, \psi \rangle \not =0$ if and only if their annihilators $E_{\varphi}$ and $E_{\psi}$ satisfy 
$E_{\varphi} \cap E_{\psi} =\{ 0 \}$. 
\end{fact}
For the proof, see I\hspace{-.1em}I\hspace{-.1em}I.2.4 in \cite{Che}. 

\begin{ex}
For a symplectic form $\omega$ on $V$, we put 
\begin{eqnarray*}
\varphi_{\omega} = \exp{\sqrt{-1}\omega}. 
%\varphi_{\omega} = \exp{\sqrt{-1}\omega} = 1 + \sqrt{-1}\omega + \frac{1}{2!}(\sqrt{-1}\omega)^{2} + \cdots + \frac{1}{m!}(\sqrt{-1}\omega)^{m} \in \wedge^{ev}V^{*} \otimes \mathbb{C}. 
\end{eqnarray*}
Then we have $E_{\varphi_{\omega}} = \{ X - \sqrt{-1}\iota_{X} \omega \big|\ X \in V \otimes \mathbb{C} \}$ and $\dim E_{\varphi_{\omega}} = n$. 
Since $\omega$ is nondegenerate, we have $\langle \varphi_{\omega}, \bar \varphi_{\omega} \rangle = \frac{(-2\sqrt{-1})^m}{m!}\omega^m \not = 0$.
Hence $\varphi_{\omega}$ is a generalized Calabi-Yau structure on $V$. The type of $\varphi_{\omega}$ is equal to $0$. 
\end{ex}

\begin{ex}
If $V$ has a complex structure $J$, then for the $\sqrt{-1}$-eigenspace $V^{1,0}$ of 
$J^{*}:V^{*} \otimes \mathbb{C} \longrightarrow V^{*} \otimes \mathbb{C}$, $\wedge^{m}V^{1,0}$ is 
one-dimensional complex vector space. 
Let $\Omega$ be a non-zero vector in $\wedge^{m}V^{1,0}$. Then, we have $E_{\Omega} = V_{0,1} \oplus V^{1,0}$ and 
$\langle \Omega, \bar \Omega \rangle = (-1)^m \Omega \wedge \bar \Omega \not = 0$. So $\Omega$ is a generalized Calabi-Yau 
structure on $V$. The type of $\Omega$ is equal to $m$. 
\end{ex}

\begin{ex}
Let $\varphi$ be a generalized Calabi-Yau structure on $V$. For each $B \in \wedge^{2}V^{*}$, 
the previous example shows that $\exp(B)\varphi$ is pure. Moreover, the bilinear form gives 
$\langle \exp(B)\varphi,\ \overline{\exp(B) \varphi} \rangle = \langle \varphi, \bar \varphi \rangle \not =0$. 
Hence $\exp(B)\varphi$ is also a generalized Calabi-Yau structure on $V$. 
The type of $\exp(B)\varphi$ coincides with that of $\varphi$. 

%Let $\varphi$ be a generalized Calabi-Yau structure on $V$. For each $B \in \wedge^{2}V^{*}$, we put 
%\begin{eqnarray*}
%\exp (B)\varphi = (1 + B + \frac{1}{2}B^{2} + \cdots ) \wedge \varphi.
%\end{eqnarray*}
%Then the bilinear form gives 
%$\langle \exp(B)\varphi,\ \overline{\exp(B) \varphi} \rangle = \langle \varphi, \bar \varphi \rangle \not =0$. 
%Moreover, since $E_{\exp(B)\varphi} = \{ X + \alpha  + \iota_{X}B \big|\ X + \alpha \in E_{\varphi} \}$, $\exp(B)\varphi$ is pure. 
%So $\exp(B)\varphi$ is also a generalized Calabi-Yau structure on $V$. 
\end{ex}

\begin{ex}
If $\varphi_{1}$ and $\varphi_{2}$ are two generalized Calabi-Yau structures on two vector spaces $V_{1}$ and $V_{2}$, and $p_{1}$, $p_{2}$ are the 
projections from the direct sum $V_{1} \oplus V_{2}$. Then $\varphi = p_{1}^{*}\varphi_{1} \wedge p_{2}^{*}\varphi_{2}$ is 
a generalized Calabi-Yau structure on the product. The type of $\varphi$ is equal to the sum type($\varphi_{1}$) + type($\varphi_{2}$). 
\end{ex}

\subsection{Generalized Calabi-Yau structures on a manifold}
Let $M$ be a smooth manifold of dimension $2n$, and consider the direct sum $TM \oplus T^{*}M$ of the tangent bundle and the cotangent bundle. 
Then there is an indefinite metric on the vector bundle $TM \oplus T^*M$ defined by 
$(X + \alpha, Y + \beta) = \frac{1}{2}( \beta (X) + \alpha (Y))$. 
\begin{defn}[\cite{Hi}]
A generalized Calabi-Yau structure on a manifold $M$ is a closed differential form $\varphi \in \Omega^{ev/od} \otimes \mathbb{C}$ which satisfies 
the following conditions. 
\begin{itemize}
\item For each $p \in M$, $\varphi_p$ is a complex pure spinor on $(T_pM \oplus T_p^*M) \otimes \mathbb{C}$. \\
\item At each point, $\langle \varphi , \bar \varphi \rangle \not = 0$. 
\end{itemize}
\end{defn}

\begin{rem}
Generalized Calabi-Yau structures were defined by Hitchin in \cite{Hi}. If a generalized Calabi-Yau structure $\varphi$ is given, 
then the annihilator $E_{\varphi}$ defines a generalized complex structure 
in the sense of Hitchin \cite{Hi}. This shows that a generalized Calabi-Yau manifold is a special case of a generalized complex manifold. 
For the detail, see Proposition 1 in \cite{Hi}.
\end{rem}

\begin{ex}
Let $M$ be an $2n$-dimensional symplectic manifold with the symplectic form $\omega$, and put 
\begin{eqnarray*}
\varphi_{\omega} = \exp {\sqrt{-1} \omega}. 
\end{eqnarray*}
Then we have $E_{\varphi_{\omega}} = \{ X - \sqrt{-1}\iota_{X} \omega \big|\ X \in T \otimes \mathbb{C} \}$ and 
$\langle \varphi_{\omega}, \bar \varphi_{\omega} \rangle = \frac{(-2\sqrt{-1})^n}{n!}\omega^n \not = 0$. 
Since $\omega$ is closed, $\varphi_{\omega}$ is also closed. Hence $\varphi_{\omega}$ is a generalized Calabi-Yau structure on $M$. 
\end{ex}

\begin{ex}
Let $M$ be a $n$-dimensional complex manifold with a non-vanishing holomorphic $n$-form $\Omega$. 
Then $\Omega$ is pure since $E_{\Omega} = T_{0,1} \oplus T^{1,0}$. In addition, the bilinear form gives 
$\langle \Omega, \bar \Omega \rangle = (-1)^n \Omega \wedge \bar \Omega$, which is non-vanishing. 
Since $\Omega$ is closed, $\Omega$ is a generalized Calabi-Yau structure on $M$. 
\end{ex}

\begin{ex}
If $B$ is a closed $2$-form on a generalized Calabi-Yau manifold $(M, \varphi)$, then $\exp (B)\varphi$
is also a closed form. By the previous example, $\exp (B)\varphi$ is pure and $\langle \exp(B)\varphi, \overline{\exp(B) \varphi} \rangle \not = 0$ 
at each point. So $\exp (B)\varphi$ is also a generalized Calabi-Yau structure on $M$. This is called the $B$-field transform of $\varphi$. 
\end{ex}

\begin{ex}
If $(M_{1}, \varphi_{1})$ and $(M_{2}, \varphi_{2})$ are two generalized Calabi-Yau manifolds and $p_{1}$, $p_{2}$ are the 
projections from the product manifold $M_{1} \times M_{2}$. Then $\varphi = p_{1}^{*}\varphi_{1} \wedge p_{2}^{*}\varphi_{2}$ is 
a generalized Calabi-Yau structure on the product. In particular, a product manifold of generalized Calabi-Yau manifolds is also a generalized Calabi-Yau manifold. 
\end{ex}

The local expression of a generalized Calabi-Yau structure is given by the following proposition by Gualtieri \cite{Gua}. 
This helps us to prove the Duistermaat-Heckman formula later. 

\begin{fact}[\cite{Gua}]\label{Gua2}
An element of a generalized Calabi-Yau manifold $(M, \varphi)$ is said to be regular if it has a neighborhood where the type of $\varphi$ is constant. 
If $p \in M$ is regular, then for sufficiently small neighborhood $U_{p}$ of $p$, there exists a complex 2-form $B + \sqrt{-1}\omega \in \Omega^{2}(U_{p}) \otimes \mathbb{C}$ such that 
\begin{eqnarray*}
\varphi = \exp (B + \sqrt{-1}\omega )\varphi_{k} \ on\ U_{p}, \\ 
\end{eqnarray*}
where $k$ is the type of $\varphi_{p}$. Moreover, $\varphi_{k}$ can be written 
\[
\varphi_{k} = \theta^{1} \wedge \cdots \wedge \theta^{k} 
\]
by some 1-forms $\theta^{1}, \cdots, \theta^{k} \in \Omega^{1}(U_{p}) \otimes \mathbb{C}$. 
\end{fact}
%%%%%%%%%%%%%%%%%%%%%%%%%%%%%%%%%%%%%%%%%%%%%%%%%%%%%%%%%%%%%%%%%%%%%%%%%%%%%%%%%%%%%%%%%%%%%%%%%%%%%%%%%%%%%%%%%%%%%%%%%%%%%
%%%%%%%%%%%%%%%%%%%%%%%%%%%%%%%%%%%%%%%%%%%%%%%%%%%%%%%%%%%%%%%%%%%%%%%%%%%%%%%%%%%%%%%%%%%%%%%%%%%%%%%%%%%%%%%%%%%%%%%%%%%%%
\section{Reduction of generalized Calabi-Yau structures}
\subsection{Generalized moment maps}
In this section we define the notion of generalized moment maps for a compact Lie group action on a generalized Calabi-Yau manifold, and construct a generalized Calabi-Yau 
structure on the reduced space. The definition of generalized moment maps for generalized complex cases is given by Lin and Tolman \cite{Lin1}. 

\begin{defn}
Let a compact Lie group $G$ with its Lie algebra $\frak{g}$ act on a generalized Calabi-Yau manifold $(M, \varphi)$ preserving $\varphi$. 
A generalized moment map is a smooth function $\mu : M \longrightarrow \frak{g}^*$ which satisfies 
\begin{itemize}
\item $\mu$ is $G$-equivariant, and 
\item $\xi_{M} - \sqrt{-1}d\mu^{\xi}$ lies in $E_{\varphi}$ for all $\xi \in \frak{g}$, where $\xi_{M}$ denotes the induced vector field on $M$ and 
$\mu^{\xi}$ is the smooth function defined by $\mu^{\xi}(p) = \mu (p)(\xi)$.
\end{itemize}
A $G$-action which preserves the generalized Calabi-Yau structure $\varphi$ is called Hamiltonian if a generalized moment map exists. 
\end{defn}

Here are some examples of generalized moment maps. 

\begin{ex}
Let $G$ act on a symplectic manifold $(M, \omega)$ preserving $\omega$, and $\mu : M \longrightarrow \frak{g}^*$ be a moment map. 
Then $G$ also preserves the generalized Calabi-Yau structure $\varphi_{\omega} = \exp \sqrt{-1}\omega$, and $\mu$ is also a generalized moment map. 
\end{ex}

\begin{ex}
Let $G$ act on a connected Calabi-Yau $n$-fold $(M, \Omega)$, where $\Omega$ is a non-vanishing holomorphic $n$-form. If the $G$-action is Hamiltonian, 
$\xi_M$ must be anti-holomorphic for all $\xi \in \frak{g}$. However induced vector fields must be real, so we have $\xi_M = 0$.
In particular, the $G$-action is trivial and the generalized moment map is regarded as a linear functional on the Lie algebra $\frak{g}$. 
\end{ex}

\begin{ex}
If $G$ acts on two generalized Calabi-Yau manifolds $(M_1, \varphi_1)$ and $(M_2, \varphi_2)$, preserving both $\varphi_1$ and $\varphi_2$. 
Let $\mu_1$ and $\mu_2$ are generalized moment maps for these actions. Then the diagonal action of $G$ on the product manifold $M_1 \times M_2$ preserves the generalized 
Calabi-Yau structure $\varphi = p_{1}^*\varphi_1 \wedge p_{2}^*\varphi_2$, where $p_{1}$ and $p_{2}$ are the 
projections from the product $M_1 \times M_2$. Moreover $\mu = \mu_1\circ p_{1} + \mu_2 \circ p_{2}$ is a generalized moment map for this action. 
\end{ex}

\subsection{Generalized Calabi-Yau structure on the reduced space}
Let a compact Lie group $G$ act on a generalized Calabi-Yau manifold $(M, \varphi)$ in a Hamiltonian way with a generalized moment map 
$\mu : M \longrightarrow \frak{g}^*.$ Suppose that $G$ acts freely on $\mu^{-1}(0).$ Then $0$ is a regular value and the quotient space
\begin{eqnarray*}
M_0 = \mu^{-1}(0)/G
\end{eqnarray*}
is a manifold. The purpose of 3.2 is to prove $\bold{Theorem\ A}$ in Introduction. By restricting to an appropriate neighborhood of $\mu^{-1}(0)$, we may assume that 
$G$ acts freely on $M$. The following lemmas are required for the proof of the theorem. 

\begin{lem}\label{dm}
Under the assumptions above, let $\frak{g}_{M}$ be the subbundle of $TM$ generated by the fundamental vector fields $\xi_{M}$ for $\xi \in \frak{g}$, 
and $d\mu$ be the subbundle of $T^{*}M$ generated by the differential $d\mu^{\eta}$ for $\eta \in \frak{g}$. 
Then we have 
\begin{description}
\item[(1)] $T_p\mu^{-1}(0) = (d\mu)_p ^{0}$, \\
\item[(2)] $\ker ({p_{0}}_*)_p = (\frak{g}_{M})_p$, and \\
\item[(3)] $T_{[p]}M_0 \cong T_p \mu^{-1}(0) / (\frak{g}_M)_{p} = (d\mu)_p ^{0} / (\frak{g}_{M})_p$, 
\end{description}
where $p \in \mu^{-1}(0)$ and $(d\mu)_p ^{0} = \left\{ X \in T_p M \big|\ (d\mu^{\xi})_p (X) = 0 \ (\xi \in \frak{g})  \right\}$ is the annihilator of $(d\mu)_p$. 
\end{lem}

\begin{proof}
For each $\xi \in \frak{g}$, the smooth function $\mu^{\xi}$ vanishes on $\mu^{-1}(0)$. So $(d\mu^{\xi})(X) = 0$ 
for all $X \in T_{p}\mu^{-1}(0)$. This implies that $T_p \mu^{-1}(0) \subset (d\mu)_p ^{0}$. 
In addition, because $\dim T_p \mu^{-1}(0) = \dim (d\mu)_p ^{0}$, the first claim holds. Since $(\frak{g}_{M})_p \subset \ker ({p_{0}}_*)_p$ and 
${p_{0}}: \mu^{-1}(0) \longrightarrow M_{0}$ is a submersion, the second claim holds. Now it is easy to see the last claim. 
\end{proof}

The following lemma will help us to prove that the reduced form does not vanish anywhere. 

\begin{lem}\label{dm2}
Under the assumptions above, let $\pi : (TM \oplus T^{*}M) \otimes \mathbb{C} \longrightarrow TM \otimes \mathbb{C}$ be the natural 
projection. Then we have 
\begin{eqnarray*}
\dim_{\mathbb{C}} (T_{p} \mu^{-1}(0) \otimes \mathbb{C}) \cap \pi ( E_{\varphi})_{p} = \dim_{\mathbb{C}} \pi (E_{\varphi})_{p} - \dim G 
\end{eqnarray*}
for each $p \in \mu^{-1}(0)$.
\end{lem}

\begin{proof}
For a subspace $W \subset (TM \oplus T^*M)_{p} \otimes \mathbb{C}$, we denote by $W^{\perp}$ 
the annihilator of $W$ with respect to the natural metric on $(TM \oplus T^*M)_{p} \otimes \mathbb{C}$. 
Then, since $E_{\varphi}$ is maximal isotropic, we have 
\begin{eqnarray*}
E_{\varphi} = E_{\varphi}^{\perp},\ {\rm and}\ W^{\perp} \cap (E_{\varphi})_{p} = (W + (E_{\varphi})_{p})^{\perp}.  
\end{eqnarray*}
If $X \in (T_{p} \mu^{-1}(0) \otimes \mathbb{C}) \cap \pi ( E_{\varphi})_{p}$, then it satisfies that 
\[
X \in \pi ( E_{\varphi})_{p},\ {\rm and}\ d\mu^{\xi}(X) = 0 
\]
for each $\xi \in \mathfrak{g}$. Thus we have $X \in \pi ((\frak{g}_{M} \otimes \mathbb{C})^{\perp} \cap E_{\varphi})_{p}$. 
Conversely, if $X \in \pi ((\frak{g}_{M} \otimes \mathbb{C})^{\perp} \cap E_{\varphi})_{p}$, then we also have $d\mu^{\xi}(X) = 0$ 
for each $\xi \in \mathfrak{g}$. So we have $X \in (T_{p} \mu^{-1}(0) \otimes \mathbb{C}) \cap \pi ( E_{\varphi})_{p}$. 
This shows that 
\[
(T_{p} \mu^{-1}(0) \otimes \mathbb{C}) \cap \pi ( E_{\varphi})_{p} = \pi ((\frak{g}_{M} \otimes \mathbb{C})^{\perp} \cap E_{\varphi})_{p}. 
\]
Since the kernel of $\pi : (TM \oplus T^{*}M)_{p} \otimes \mathbb{C} \longrightarrow T_{p}M \otimes \mathbb{C}$ 
is equal to $T_{p}^{*}M \otimes \mathbb{C}$, we have 
\begin{eqnarray*}
(\frak{g}_{M} \otimes \mathbb{C})_{p}^{\perp} \cap (E_{\varphi})_{p} \cap T_{p}^*M \otimes \mathbb{C} &=& ((\mathfrak{g}_{M} \otimes \mathbb{C}) + E_{\varphi})_{p}^{\perp} \cap T_{p}^{*}M \otimes \mathbb{C} \\
                                                                                                      &=& \pi((\mathfrak{g}_{M} \otimes \mathbb{C}) + E_{\varphi})_{p}^{0} \\
                                                                                                      &=& \pi(E_{\varphi})_{p}^{0}, 
\end{eqnarray*}
and thus
\begin{eqnarray*}
\dim_{\mathbb{C}} \pi ( (\frak{g}_{M} \otimes \mathbb{C})^{\perp} \cap E_{\varphi})_{p} = \dim_{\mathbb{C}} (\frak{g}_{M} \otimes \mathbb{C})_{p}^{\perp} \cap (E_{\varphi})_{p} - \dim_{\mathbb{C}} \pi(E_{\varphi})_{p}^{0}. 
\end{eqnarray*}
In addition, by $(\frak{g}_{M} \otimes \mathbb{C})_{p} \cap (E_{\varphi})_{p} = \{ 0 \}$, we obtain the dimension 
\begin{eqnarray*}
\dim_{\mathbb{C}} (\frak{g}_{M} \otimes \mathbb{C})_{p}^{\perp} \cap (E_{\varphi})_{p} &=& \dim_{\mathbb{C}} ((\frak{g}_{M} \otimes \mathbb{C}) + E_{\varphi})_{p}^{\perp} \\
                                                                                       &=& \dim M - \dim G.
\end{eqnarray*}
Hence we have 
\begin{eqnarray*}
\dim_{\mathbb{C}} \left (T_{p} \mu^{-1}(0) \otimes \mathbb{C} \right) \cap \pi ( E_{\varphi})_{p} &=& \dim_{\mathbb{C}}\pi ((\frak{g}_{M} \otimes \mathbb{C})^{\perp} \cap (E_{\varphi}))_{p} \\
                                                                                                  &=& \dim_{\mathbb{C}} (\frak{g}_{M} \otimes \mathbb{C})_{p}^{\perp} \cap (E_{\varphi})_{p} - \dim_{\mathbb{C}} \pi(E_{\varphi})_{p}^{0} \\
                                                                                                  &=& \dim_{\mathbb{C}} \pi (E_{\varphi})_{p} - \dim G, 
\end{eqnarray*}
this completes the proof. 
\end{proof}

\begin{proof}[Proof of Theorem \ref{main1}]
For each $p \in \mu^{-1}(0)$, we denote by $(\varphi_{s})_{p}$ the $s$-th degree component of $\varphi_{p} \in \wedge^{ev/od}T_{p}^{*}M \otimes \mathbb{C}$. 
Then, by the definition of the generalized moment map, we have 
\[
\iota_{\xi_{M}}\varphi_{s} - \sqrt{-1}d\mu^{\xi} \wedge \varphi_{s-2} = 0 
\]
for each $\xi \in \mathfrak{g}$. Moreover, the identity $T_{p}\mu^{-1}(0) = (d\mu)_{p}^{0}$ in Lemma \ref{dm} tells us that the $(s-1)$-form 
$\iota_{(\xi_{M})_{p}}(\varphi_{s})_{p}$ vanishes on $T_{p}\mu^{-1}(0)$. So by identifying the tangent space $T_{[p]}M_0$ with 
$T_p \mu^{-1}(0) / (\frak{g}_M)_{p}$ (see Lemma \ref{dm}, (3)), we obtain a well-defined complex $s$-form
$(\tilde \varphi_{s})_{[p]}$ on $T_{[p]}M_{0}$ by 
\[
(\tilde \varphi_{s})_{[p]}([X_{1}], \cdots, [X_{s}]) = (i_{0}^{*}\varphi_{s})_{p}(X_{1}, \cdots, X_{s}), 
\]
where $X_{1}, \cdots, X_{s} \in T_{p}\mu^{-1}(0)$. Thus we have a complex form 
$(\tilde \varphi)_{[p]} \in \wedge^{ev/od}T_{[p]}^{*}M_{0} \otimes \mathbb{C}$ defined by 
\[
(\tilde \varphi)_{[p]} = (\tilde \varphi_{k})_{[p]} + (\tilde \varphi_{k+2})_{[p]} + \cdots, 
\]
where $k$ is the type of $\varphi_{p}$. $G$-invariance of the form $\varphi$ tells us that the definition of $(\widetilde \varphi)_{[p]}$ 
does not depend on a representative $p \in \mu^{-1}(0)$. 
So we get the reduced form $\widetilde \varphi \in \Omega^{ev/od} \otimes \mathbb{C}$. 
It is clear that $\widetilde \varphi$ satisfies thet $p_{0}^* \widetilde \varphi = i_{0}^{*} \varphi$ and $d\widetilde \varphi = 0$. 

Next we shall show that $(\widetilde \varphi)_{[p]} \not = 0$. It is sufficient to show that $(i_{0}^{*} \varphi_{k})_{p} \not = 0$. 
Suppose that $\dim M = 2n$ and $\dim G = l$. Then Lemma \ref{dm2} tells us 
\begin{eqnarray*}
\dim_{\mathbb{C}} (T_{p} \mu^{-1}(0) \otimes \mathbb{C}) \cap \pi ( E_{\varphi})_{p}  = 2n - k - l. 
\end{eqnarray*}
So we can take a basis 
\begin{eqnarray*}
e_{1}, \cdots, e_{2n-k-l}, u_{1}, \cdots, u_{k}, v_{1},\cdots, v_{l} 
\end{eqnarray*}
of $T_{p} M \otimes \mathbb{C}$, where $\{e_{1}, \cdots, e_{2n-k-l}, u_{1}, \cdots, u_{k}\}$ is a basis of $T_{p}\mu^{-1}(0) \otimes \mathbb{C}$, 
and $\{e_{1}, \cdots, e_{2n-k-l}, v_{1},\cdots, v_{l}\}$ is a basis of $\pi (E_{\varphi})$. 
Since $(\varphi_{k})_{p} \not = 0$, so we have 
\begin{eqnarray*}
(\varphi_{k})_{p}(u_{1}, \cdots , u_{k}) \not = 0. 
\end{eqnarray*}
This shows that $(i_{0}^{*} \varphi_{k})_{p} \not = 0$. 

Now we say that an element $\tilde X + \tilde \alpha \in (TM_{0} \oplus T^{*}M_{0})_{[p]} \otimes \mathbb{C}$ 
satisfies the compatiblity condition if there exists $X \in T_{p}\mu^{-1}(0) \otimes \mathbb{C}$ and 
$\alpha \in T_{p}^*M \otimes \mathbb{C}$ such that 
$({p_{0}}_{*})_{p}X = \tilde X$, 
$p_{0}^{*}\tilde \alpha = i_{0}^{*}\alpha$, and that 
$({i_{0}}_{*})_{p}(X) + \alpha \in (E_{\varphi})_{p}$. 
We denote by $E_{0}$ the set of elements 
$\tilde X + \tilde \alpha \in (TM_{0} \oplus T^{*}M_{0})_{[p]} \otimes \mathbb{C}$ which satisfy the compatiblity condition. 
Then, for each $\tilde X + \tilde \alpha \in E_{0}$, we have 
\begin{eqnarray*}
p_{0}^{*}(\iota_{\tilde X}\tilde \varphi + \tilde \alpha \wedge \tilde\varphi) &=& i_{0}^{*}(\iota_{({i_{0}}_{*})X} \varphi + \alpha \wedge \varphi) \\
                                                                        &=& 0. 
\end{eqnarray*}
So we can see $E_{0} \subset E_{\tilde \varphi}$ because $p_{0}$ is a submersion. Moreover, since $E_{\tilde \varphi}$ is isotropic, we have $\dim_{\mathbb{C}} E_{0} \leq \dim_{\mathbb{C}} E_{\tilde \varphi} \leq 2(n-l)$. 
Let us show the equality $\dim_{\mathbb{C}} E_{0} = 2(n-l)$. Since 
$\dim_{\mathbb{C}} (T_{p} \mu^{-1}(0) \otimes \mathbb{C}) \cap \pi ( E_{\varphi})_{p}  = 2n - k - l$, we can take 
\[
X_{1} + \alpha_{1}, \cdots, X_{2n-l-k} + \alpha_{2n-l-k} \in E_{\varphi}, 
\]
which are linearly independent and $X_{i} \in T_{p}\mu^{-1}(0) \cap \pi(E_{\varphi})_{p}$ for $i = 1, \cdots, 2n-l-k$. 
Since 
\[
\iota_{\xi_{M}}\alpha_{i} = (\alpha_{i}, \xi_{M}) = (X_{i} + \alpha_{i}, \xi_{M} - d\mu^{\xi}) = 0 
\]
for each $\xi \in \mathfrak{g}$, $\alpha_{i}$ descends to a form $\tilde \alpha_{i} \in \wedge^{ev/od}T_{[p]}^{*}M_{0}\otimes \mathbb{C}$. 
If we take 
\[
\tilde X_{i} = ({p_{0}}_{*})_{p}X, 
\]
then we have $\tilde X_{i} + \tilde \alpha_{i} \in E_{0}$. Furthermore, since $\ker({p_{0}}_{*})_{p} = (\mathfrak{g}_{M})_{p}$ 
has dimension $l$, and it is contained in $T_{p}\mu^{-1}(0) \cap \pi(E_{\varphi})_{p}$, 
so we may assume that 
\[
\tilde X_{1} + \tilde \alpha_{1}, \cdots, \tilde X_{2(n-l)-k} + \tilde \alpha_{2(n-l)-k} 
\]
are linearly independent. 

On the other hand, by Fact \ref{fact:Gua}, 
we can take $\theta^{1}, \cdots, \theta^{k} \in T_{p}^{*}M \otimes \mathbb{C}$ which 
satisfy 
\[
(\varphi_{k})_{p} = \theta^{1} \wedge \cdots \wedge \theta^{k}. 
\]
Then, since $(\varphi_{k})_{p}$ satisfies $\iota_{\xi_{M}} (\varphi_{k})_{p} = 0$ for each $\xi \in \mathfrak{g}$, 
so does $\theta^{i}$ for $i = 1, \cdots, k$. Hence $\theta^{i}$ descends to a $1$-form $\tilde \theta^{i} \in \wedge^{ev/od}T_{[p]}^{*}M_{0}\otimes \mathbb{C}$. 
Then $\tilde \theta^{i} \in E_{0}$, and 
\begin{eqnarray*}
p_{0}^{*}(\tilde \theta^{1} \wedge \cdots \wedge \tilde \theta^{k}) &=& i_{0}^{*}(\theta^{1} \wedge \cdots \wedge \theta^{k}) \\
                                                                    &=& i_{0}^{*}((\varphi_{k})_{p}) \\
                                                                    &\not=& 0. 
\end{eqnarray*}
This shows that $\tilde \theta^{1}, \cdots, \tilde \theta^{k}$ are linearly independent. 
Thus we have 
\[
\dim_{\mathbb{C}} E_{0} = 2(n-l),\ {\rm and}\ E_{0} = E_{\tilde\varphi}, 
\]
in particular $E_{\tilde \varphi}$ is maximal isotropic. 

Furthermore, since $E_{\varphi}$ does not have a real vector except for $0$, neither does $E_{0}$. 
So we also have 
\[
(E_{\widetilde \varphi})_{[p]} \cap (\bar{E}_{\widetilde \varphi})_{[p]} = \{ 0 \}. 
\]
This shows that $\widetilde \varphi$ is a generalized Calabi-Yau structure on $M_0$. 

The last claim is clear because ${\rm type}(\varphi_{p}) = {\rm type}((\widetilde{\varphi})_{[p]}) = k$.  
\end{proof}

\begin{rem}
The reduction for other levels can be done by taking the coadjoint orbit. The detailed statement is as follows. 
Let a compact Lie group $G$ act on a generalized Calabi-Yau manifold $(M, \varphi)$ Hamiltonian way with a generalized moment map 
$\mu : M \longrightarrow \frak{g}^*$. For each $a \in \frak{g}^{*}$,  $\mathcal{O}_{a}$ denotes 
the coadjoint orbit of $a$. Suppose that $G$ acts on $\mu^{-1}(\mathcal{O}_{a})$ freely. Then the quotient space $M_{a} = \mu^{-1}(\mathcal{O}_{a})/G$ is a manifold 
and has unique generalized Calabi-Yau structure $\widetilde \varphi$ which satisfies that 
\begin{eqnarray*}
p_{a}^* \widetilde \varphi = i_{a}^{*} \varphi 
\end{eqnarray*}
and
\begin{eqnarray*}
{\rm type}(\varphi_p) = {\rm type}(\widetilde \varphi_{[p]})
\end{eqnarray*}
for all $p \in \mu^{-1}(\mathcal{O}_{a})$, where $i_{a} : \mu^{-1}(\mathcal{O}_{a}) \longrightarrow M$ is 
the inclusion and $p_{a}: \mu^{-1}(\mathcal{O}_{a}) \longrightarrow M_{a}$ is the natural projection. 
In addition, we have $\dim M_{a} = \dim M + \dim \mathcal{O}_{a} - 2\dim G$. 
\end{rem}

\begin{ex}
Let $G$ act on a symplectic manifold $(M, \omega)$ preserving $\omega$, and let $\mu : M \longrightarrow \frak{g}^*$ be a moment map. 
Then $G$ also acts on $(M, \varphi_{\omega})$ Hamiltonian way and $\mu$ is a generalized moment map. Moreover if we assume that $G$ acts freely 
on $\mu^{-1}(0)$, then we get the reduced symplectic structure $\widetilde \omega$ and the reduced generalized Calabi-Yau structure $\widetilde \varphi_{\omega}$
on the reduced space $M_0$. Then $\widetilde \varphi_{\omega}$ coincides with the generalized Calabi-Yau structure 
$\varphi_{\widetilde \omega}$ induced by the reduced symplectic structure $\widetilde \omega$. 
\end{ex}

\begin{ex}
Let $G$ act on a Calabi-Yau manifold $(M, \Omega)$. If the $G$-action is Hamiltonian, 
then the action is trivial and the generalized moment map $\mu$ is regarded as a linear functional on the Lie algebra $\frak{g}$. 
So the reduced space $M_{0}$ coincides with either $M$ or the empty set. 
\end{ex}

\begin{rem}
Lin and Tolman showed the existence of a generalized complex structure on the reduced space in \cite{Lin1}. The generalized complex structure induced by the reduced 
generalized Calabi-Yau structure coincides with the reduced generalized complex structure from the generalized complex structure induced by the original generalized Calabi-Yau 
structure. 
\end{rem}

\subsection{Relationship to Bergman kernels}
We introduce a Hamiltonian action on a generalized Calabi-Yau structure which is not induced from either a symplectic structure or a Calabi-Yau one here. 
Let $D \subset \mathbb{C}^{m+n}$ be a Reinhardt bounded domain, that is, a bounded domain which the standard action of $(m+n)$-dimensional torus $T^{m+n}$ on $\mathbb{C}^{m+n}$ 
leaves $D$ invariant. For each $w = (w_{1}, \cdots , w_{m}) \in \mathbb{C}^{m}$, $D_{w}$ denotes the slice of $D$ at $w$, 
\[
D_{w} = \left\{ (z_{1},\cdots , z_{m+n}) \in D\ \big|\ z_{j} = w_{j}\quad (j = 1, \cdots , m) \right\}.
\]
If the slice $D_{w}$ is not empty, we can regard $D_{w}$ as a Reinhardt bounded domain in $\mathbb{C}^{n}$ naturally. Let 
\[
K_{w}(z)=K_{w}(z, z) : D_{w} \longrightarrow \mathbb{R} 
\]
be the Bergman kernel function of $D_{w}$, and $\Omega_{w} = \frac{\sqrt{-1}}{2}\partial \bar \partial \log K_{w}$ be the K${\rm {\ddot a}}$hler form 
of the Bergman metric on $D_{w}$. Then the natural action of $S^{1}$ on $D_{w}$ preserves $\Omega_{w}$, and 
\[
\mu_{w} = -\frac{1}{2}\sum_{j = 1}^{n} z_{j} \frac{\partial}{\partial z_{j}}\big( \log K_{w} \big) 
\]
is a moment map for this action. Note that the function $\mu_{w}$ is real and $S^{1}$-invariant since the real function $\log K_{w}$ is $S^{1}$-invariant and 
the fundamental vector field $\xi$ induced by the $S^{1}$-action is given by 
\[
\xi = \sqrt{-1}\sum_{j=1}^{n}\left\{ z_{j}\frac{\partial}{\partial z_{j}} - \bar z_{j}\frac{\partial}{\partial \bar z_{j}} \right\}. 
\]

Now we assume that the Bergman kernel $K_{w}$ depends smoothly on $w$. Then we can define a smooth function $K$ on $D$ by 
\[
K(w, z) = K_{w}(z) : D \longrightarrow \mathbb{R}, 
\]
and a complex form $\varphi$ on $D$ by 
\[
\varphi = dw_{1}\wedge \cdots \wedge dw_{m} \wedge \exp \sqrt{-1}\Omega, 
\]
where $\Omega = \frac{\sqrt{-1}}{2}\partial \bar \partial \log K$. 
It is easy to see that the complex form $\varphi$ is a generalized Calabi-Yau structure on $D$, and the $S^{1}$-action on $D$ defined by 
\[
e^{\sqrt{-1}\theta}(w_{1}, \cdots , w_{m}, z_{1}, \cdots , z_{n}) = (w_{1},\cdots , w_{m}, e^{\sqrt{-1}\theta}z_{1}, \cdots , e^{\sqrt{-1}\theta}z_{n}), 
\]
preserves $\varphi$. 

\begin{thm}
Let $\mu$ be a smooth function on $D$ defined by 
\[
\mu(w_{1}, \cdots, w_{m},z_{1}, \cdots, z_{n}) = \mu_{w}(z_{1}, \cdots , z_{n}) = -\frac{1}{2}\sum_{j = 1}^{n} z_{j} \frac{\partial}{\partial z_{j}}\big( \log K \big). 
\]
Then the function $\mu$ is a generalized moment map for the $S^{1}$ action on $D$ defined above. 
\end{thm}

\begin{proof}
Let $\xi$ be the fundamental vector field for this action. Then $S^{1}$-invariance of the function $\log K$ implies that $\mu$ is a 
$S^{1}$-invariant real-valued function. By simple calculation, we have 
\begin{eqnarray*}
\iota_{\xi}\Omega \big( \frac{\partial}{\partial z_{i}} \big) &=& \Omega \big( \sqrt{-1}\sum_{j=1}^{n} \left\{ z_{j}\frac{\partial}{\partial z_{j}} - \bar z_{j} \frac{\partial}{\partial \bar z_{j}} \right\}, \frac{\partial}{\partial z_{i}} \big) \\
                                                              &=& \sqrt{-1}\sum_{j=1}^{n} \bar z_{j} \left( \frac{\sqrt{-1}}{2}\frac{\partial^{2}}{\partial z_{i} \partial \bar z_{j}} \left( \log K \right) \right) \\
                                                              &=& \frac{\partial}{\partial z_{i}} \left( -\frac{1}{2} \sum_{j=1}^{n} \bar z_{j} \frac{\partial}{\partial \bar z_{j}} \left(\log K \right) \right) \\ 
                                                              &=& \frac{\partial}{\partial z_{i}} \left( -\frac{1}{2} \sum_{j=1}^{n} z_{j} \frac{\partial}{\partial z_{j}} \left(\log K \right) \right) \\
                                                              &=& \frac{\partial \mu}{\partial z_{i}}, 
\end{eqnarray*}
and $\iota_{\xi}\Omega \big( \frac{\partial}{\partial w_{i}} \big) = \frac{\partial \mu}{\partial w_{i}}$ similarly. Hence we have $d\mu = \iota_{\xi}\Omega$, and 
we can check easily that $\mu$ is a generalized moment map for this action. 
\end{proof}

\begin{ex}\label{ex1}
Let $D$ be an $(m+n)$-dimensional polydisc, 
\[
D = (D^{1})^{m+m} = \{(z_{1}, \cdots, z_{m+n})\ \big|\ |z_{j}| < 1\quad (j = 1,\cdots, m+n)\}. 
\]
For each $w \in (D^{1})^{m} = \{ (w_{1}, \cdots, w_{m}) \in \mathbb{C}^{m}\ \big|\ |w_{j}| < 1\quad (j = 1, \cdots m) \}$, 
$D_{w}$ denote the slice of $D$ at $w$,
\[
D_{w} = \{ (z_{1}, \cdots, z_{n}) \in \mathbb{C}^{n}\ \big|\ |z_{j}| < 1\quad (j = 1, \cdots n) \}. 
\]
Then $D_{w}$ is a polydisc on $\mathbb{C}^{n}$, and 
\[
K_{w} = \frac{1}{\pi^{n}} \frac{1}{\Pi_{j= 1}^{n}(1 - |z_{j}|^{2})^{2}}
\]
is the Bergman Kernel function of $D_{w}$. Since the Bergman kernel $K_{w}$ does not depend on $w$, 
\[
K(w, z) = K_{w}(z) : D \longrightarrow \mathbb{R},
\]
is a smooth function on $D$, and thus we get a generalized Calabi-Yau atructure on $D$, 
\[
\varphi = dw_{1} \wedge \cdots \wedge dw_{m} \wedge \exp \sqrt{-1}\Omega, 
\]
where $\Omega = \frac{\sqrt{-1}}{2}\partial \bar \partial \log K$. The natural $S^{1}$-action defined above 
preserves $\varphi$, and we have a generalized moment map $\mu$ for this action, 
\[
\mu = -\sum_{j=1}^{n}\frac{|z_{j}|^{2}}{1 - |z_{j}|^{2}}. 
\]

On the other hand, since the total space $D$ and the parameter space $(D^{1})^{m}$ are also Reinhardt bounded domains, 
They have K$\rm{\ddot a}hler$ forms induced by their Bergman kernels. So they have also generalized Calabi-Yau 
structures induced by their K$\rm{\ddot a}hler$ forms, and They are preserved by the natural $S^{1}$-actions on them. 
By simple calculations, we get moment maps for their actions, 
\[
\mu_{D} = -\left( \sum_{ i = 1}^{m}\frac{|w_{i}|^{2}}{1-|w_{j}|^{2}} + \sum_{j=1}^{n}\frac{|z_{j}|^{2}}{1 - |z_{j}|^{2}} \right)
\]
on $D$, and 
\[
\mu_{D^{m}} = -\sum_{ i = 1}^{m}\frac{|w_{i}|^{2}}{1-|w_{j}|^{2}} 
\]
on $D^{m}$. Then they satisfy the following additive relation; 
\[
\mu_{D} = \mu_{D^{m}} + \mu. 
\]
\end{ex}

\begin{ex}\label{ex2}
Let $D$ be an $(m+n)$-dimensional complex ball 
\[
D = D^{m+n} = \{ (w_{1}, \cdots, w_{m}, z_{1},\cdots z_{n}) \in \mathbb{C}^{m+n}\ \big|\ \sum_{j=1}^{m}|w_{j}|^{2} + \sum_{j=1}^{n}|z_{j}|^{2} < 1 \}. 
\]
For each $w \in D^{m} = \{ (w_{1}, \cdots, w_{m}) \in \mathbb{C}^{m}\ \big|\ \sum_{j=1}^{m}|w_{j}|^{2} < 1 \}$, 
$D_{w}$ denote the slice of $D$ at $w$,
\[
D_{w} = \{ (z_{1}, \cdots, z_{n}) \in \mathbb{C}^{n}\ \big|\ \sum_{j=1}^{n}|z_{j}|^{2} < 1 - \sum_{j=1}^{m}|w_{j}|^{2} \}. 
\]
Then $D_{w}$ is also a complex ball on $\mathbb{C}^{n}$, and 
\[
K_{w} = \frac{n!}{\pi^{n}} \frac{1 - \sum_{j=1}^{m}|w_{j}|^{2}}{(1 - \sum_{j=1}^{m}|w_{j}|^{2} - \sum_{j=1}^{n}|z_{j}|^{2})^{n+1}}
\]
is the Bergman Kernel function of $D_{w}$. Since the Bergman kernel $K_{w}$ depends smoothly on $w$, 
\[
K(w, z) = K_{w}(z) : D \longrightarrow \mathbb{R},
\]
is a smooth function on $D$, and thus we get a generalized Calabi-Yau atructure on $D$, 
\[
\varphi = dw_{1} \wedge \cdots \wedge dw_{m} \wedge \exp \sqrt{-1}\Omega, 
\]
where $\Omega = \frac{\sqrt{-1}}{2}\partial \bar \partial \log K$. The natural $S^{1}$-action on $D$ 
preserves $\varphi$, and we have a generalized moment map $\mu$ for this action, 
\[
\mu = -\frac{n+1}{2} \frac{1 - \sum_{j = 1}^{m}|w_{j}|^{2}}{1 - ( \sum_{j=1}^{m}|w_{j}|^{2} + \sum_{j=1}^{n}|z_{j}|^{2})}. 
\]

As in the case of the previous example, we have moment maps for the natural actions of $S^{1}$ on $D$ and $D^{m}$ which are derived from 
their Bergman kernels, 
 
\[
\mu_{D} = -\frac{m+n+1}{2}\frac{1}{1-(\sum_{j=1}^{m}|w_{j}|^{2} + \sum_{j=1}^{n}|z_{j}|^{2})} 
\]
on $D$, and
\[
\mu_{D^{m}} = -\frac{m+1}{2} \frac{1}{1 - \sum_{j=1}^{m}|w_{j}|^{2}} 
\]
on $D^{m}$. They have the following multiplicative relation; 
\[
\mu_{D} = -\frac{2(m+n+1)}{(m+1)(n+1)}\mu_{D^{m}} \cdot \mu. 
\]

\end{ex}
%%%%%%%%%%%%%%%%%%%%%%%%%%%%%%%%%%%%%%%%%%%%%%%%%%%%%%%%%%%%%%%%%%%%%%%%%%%%%%%%%%%%%%%%%%%%%%%%%%%%%%%%%%%%%%%%%%%%%%%%%%%%%
%%%%%%%%%%%%%%%%%%%%%%%%%%%%%%%%%%%%%%%%%%%%%%%%%%%%%%%%%%%%%%%%%%%%%%%%%%%%%%%%%%%%%%%%%%%%%%%%%%%%%%%%%%%%%%%%%%%%%%%%%%%%%
\section{The Duistermaat-Heckman formula}
\subsection{The Duistermaat-Heckman measures and the reduced volumes}
Let $(M, \varphi)$ be a $2n$-dimensional connected generalized Calabi-Yau manifold which has constant type $k$, and suppose that compact $l$-torus $T$ acts on $M$ effectively and in a Hamiltonian way. 
In addition, we assume that the generalized moment map $\mu$ is proper. 
Then we have a natural volume form
\begin{eqnarray*}
dm = \frac{(\sqrt{-1})^n}{2^{n-k}} \langle \varphi , \bar \varphi \rangle .
\end{eqnarray*}
The volume form $dm$ defines a measure on $M$. 
Our second purpose is to prove the Duistermaat-Heckman formula in this case. 

Let $\frak{t}$ denote the Lie algebra of $T$, and $\frak{t}_{\rm reg}^{*}$ denote the subset of $\frak{t}^{*}$ 
consisting of the regular values of $\mu$. If $a \in \frak{t}^{*}$ is a regular value of $\mu$ and $p \in \mu^{-1}(a)$, 
then the stabilizer group 
\[
T_{p} = \{ g \in T \big|\ g \cdot p =p \}
\]
is finite. So if $T$-action on $\mu^{-1}(a)$ is not free, the quotient space $M_{a} = \mu^{-1}(a)/T$ is an orbifold. In this case, 
There exists a complex differential form on $M_{a}$ which, in each local representation is a generalized Calabi-Yau structure on $\mathbb{R}^{2(n-l)}$, 
and satisfies 
\[
p_{a}^{*} \widetilde \varphi = i_{a}^{*} \varphi, 
\]
where $i_{a} : \mu^{-1}(a) \longrightarrow M$ is the inclusion and $p_{a}: \mu^{-1}(a) \longrightarrow M_a$ is the natural projection. 
We call it a generalized Calabi-Yau structure on an orbifold $M_a$. 

Since $\mu$ is proper, $\frak{t}_{\rm reg}^{*}$ is a dense open subset, 
and $\frak{t}^{*} \setminus \frak{t}_{\rm reg}^{*}$ has measure 0 because of Sard's theorem. The following lemma is due to Appendix B in \cite{Gui}. 

\begin{lem}\label{null}
Suppose that $M$ is connected and $T$ acts on $M$ effectively. Then the set $M_{{\rm free}}$ on which $T$ acts freely is equal to the complement of 
a locally finite union of submanifolds of codimension $\geq$ 2. 
In particular $M_{{\rm free}}$ is open, connected, dense, and $M \setminus M_{{\rm free}}$ has measure $0$. Also $(\mu_{*})_{p}$ is surjective for all $p \in M_{{\rm free}}$. 
\end{lem}

Now we consider the normalized Haar measure $dt$ on $T$. Then the measure $dt$ induces the Lebesgue measure $dX$ on its Lie algebra $\mathfrak{t}$, and we obtain 
the dual Lebesgue measure $d\zeta$ on $\frak{t}^{*}$. 
The assumption that $\mu$ is proper implies that the pushforward $\mu_{*}(dm)$ of $dm$ under $\mu$ defines a measure in $\frak{t}^{*}$. 
We call it the Duistermaat-Heckman measure. In view of Lemma \ref{null}, we obtain $M \setminus M_{{\rm free}}$ has measure $0$ and 
$\mu|_{M_{{\rm free}}} : M_{{\rm free}} \longrightarrow \mathfrak{t}^{*}$ is a submersion. This shows that 
$\mu_{*}(dm)$ is absolutely continuous with respect to the Lebesgue measure $d\zeta$. So there exists a Borel measurable function $f$ 
on $\frak{t}^{*}$ which satisfies 
\begin{eqnarray*}
\mu_{*}(dm) = f d\zeta. 
\end{eqnarray*}
The corresponding Duistermaat-Heckman formula is stated in $\bold{Theorem\ B}$ in Introduction. 
For the proof, we need the following lemma. 

\begin{lem}\label{loc} 
For each regular point $p \in M$ of the generalized moment map $\mu$, there exists a neighborhood $U_{p}$ of $p$ and a complex $2$-form $B + \sqrt{-1}\omega \in \Omega^{2}(U_{p}) \otimes \mathbb{C}$ such that 
$\varphi = \exp(B + \sqrt{-1}\omega)\varphi_{k}$ on $U_{p}$, and $\iota_{\xi_{M}}\omega = d\mu^{\xi}$ for all $\xi \in \frak{t}$. 
\end{lem}

\begin{proof}
By Fact \ref{Gua2}, there exists a neighborhood $U_{p}$ and a complex $2$-form $\widetilde B + \sqrt{-1}\widetilde \omega \in \Omega^{2}(U_{p}) \otimes \mathbb{C}$ such that 
$\varphi = \exp(\widetilde B + \sqrt{-1}\widetilde \omega)\varphi_{k}$ on $U_{p}$. Moreover, there exists a local frame 
$\theta^{1}, \cdots , \theta^{2n}$ of $\wedge^{1}T^{*}M$ such that $\varphi_{k} = \theta^{1} \wedge \cdots \wedge \theta^{k}$ on $U_{p}$. 
So we may assume that $\widetilde B + \sqrt{-1}\widetilde \omega$ can be written 
\[
\widetilde B + \sqrt{-1}\widetilde \omega = \sum_{i,j > k} c_{ij}\theta^{i} \wedge \theta^{j}, 
\]
where $c_{ij}$ is a smooth complex function on $U_{p}$. In addition, since $p$ is a regular point, so $(\mathfrak{t}_{M})_{p}$ has dimension $l$. Hence we may 
assume that $\mathfrak{t}_{M}$ has dimension $l$ on $U_{p}$. 

Now consider the dual basis $\{ X_{1}, \cdots , X_{2n} \}$ of $\{ \theta^{1}, \cdots , \theta^{2n} \}$, 
and take an arbitrary Riemannian metric on $M$. Then we can define a complex $1$-forms $\eta^{1}, \cdots , \eta^{k}$ on $U_{p}$ defined by 
\[
\eta^{i}(\xi_{M}) = \sqrt{-1}d\mu^{\xi}(X_{i}) 
\]
for $\xi_{M} \in \frak{t}_{M}$, and vanishes on the orthogonal complement of $\mathfrak{t}_{M}$. 
Then we define a complex $2$-form $B + \sqrt{-1}\omega$ on $U_{p}$ by 
\[
B + \sqrt{-1}\omega = \widetilde B + \sqrt{-1}\widetilde \omega + \sum_{s = 1}^{k}\eta^{s} \wedge \theta^{s}. 
\]
It is clear that $\varphi = \exp(B + \sqrt{-1}\omega)\varphi_{k}$ on $U_{p}$ and 
\begin{eqnarray*}
\iota_{\xi_{M}}(B + \sqrt{-1}\omega)(X_{i}) &=& ( \sum_{i, j > k}c_{ij}\theta^{i} \wedge \theta^{j} )(\xi_{M}, X_{i}) + ( \sum_{s= 1}^{k}\eta^{s} \wedge \theta^{s} )(\xi_{M}, X_{i}) \\
                                            &=& ( \sum_{s= 1}^{k}\eta^{s} \wedge \theta^{s} )(\xi_{M}, X_{i}) \\
                                            &=& \sum_{s = 1}^{k}\eta^{s}(\xi_{M})\theta^{s}(X_{i}) \\
                                            &=& \eta^{i}(\xi_{M}) \\
                                            &=& \sqrt{-1}d\mu^{\xi}(X_{i}), 
\end{eqnarray*}
for each $\xi \in \mathfrak{t}$ and $i=1,\cdots, k$. On the other hand, since $\xi_{M} - \sqrt{-1}d\mu^{\xi} \in E_{\varphi}$ 
for each $\xi \in \frak{t}$, so we have 
\[
(\iota_{\xi_{M}}( B + \sqrt{-1}\omega ) - \sqrt{-1}d\mu^{\xi}) \wedge \varphi_{k} = 0. 
\]
Thus for $i = k+1, \cdots, 2n$, we obtain 
\begin{eqnarray*}
0 &=& \iota_{X_{i}}\left( (\iota_{\xi_{M}}( B + \sqrt{-1}\omega ) - \sqrt{-1}d\mu^{\xi}) \wedge \varphi_{k} \right) \\
   &=& (\iota_{\xi_{M}}( B + \sqrt{-1}\omega ) - \sqrt{-1}d\mu^{\xi})(X_{i})\varphi_{k}, 
\end{eqnarray*}
and  
\[
\iota_{\xi_{M}}( B + \sqrt{-1}\omega )(X_{i}) = \sqrt{-1}d\mu^{\xi}(X_{i}). 
\] 
This shows that 
\[
\iota_{\xi_{M}}(B + \sqrt{-1}\omega) = \sqrt{-1}d\mu^{\xi}, 
\]
and in particular we have $\iota_{\xi_{M}}\omega = d\mu^{\xi}$. 

\end{proof}

\begin{proof}[Proof of Theorem \ref{main2}]
Let $a \in \frak{t}_{{\rm reg}}^{*}$ be an arbitrary regular value of $\mu$ and $U$ be a convex neighborhood of $a$ contained in $\frak{t}_{{\rm reg}}^{*}$. 
Since $\frak{t}_{{\rm reg}}^{*}$ is an open set of $\frak{t}^{*}$, there exists such a neighborhood. 
Now consider a $T$-invariant connection for the fibration $\mu : \mu^{-1}(U) \longrightarrow U$. 
For each $p \in \mu^{-1}(U)$, draw the horizontal curves lying over the straight lines through $a$ and $b = \mu(p)$. This defines a $T$-equivariant projection 
$\Phi : \mu^{-1}(U) \longrightarrow \mu^{-1}(a)$ 
such that for each $b \in U$ the restriction $\Phi \big|_{\mu^{-1}(b)} : \mu^{-1}(b) \longrightarrow \mu^{-1}(a)$ is a $T$-equivariant diffeomorphism and 
\[
\mu \times \Phi : \mu^{-1}(U) \longrightarrow U \times \mu^{-1}(a) 
\]
is a trivialization. Using this trivialization and Fubini theorem, we have that $f(a)$ is equal to the 
volume of $\mu^{-1}(a)$ with respect to the quotient of $dm$ by $\mu^{*}d\zeta$. In addition, by Lemma \ref{loc} 
$dm/\mu^{*}d\zeta$ is locally given by the ($2n-l$)-form
\[
i_{a}^{*}(\varphi_{k} \wedge \bar \varphi_{k}) \wedge \frac{1}{(n-k-l)!}(i_{a}^{*}\omega)^{n-k-l} \wedge \eta, 
\]
where $\omega$ is a $2$-form given by the lemma above and $\eta$ is an $l$-form which on the $T$-orbits takes the value $\pm 1$ on an $l$-tuple $(X_{1}, \cdots , X_{l})$ 
such that $dX(X_{1}, \cdots , X_{l}) = 1$. 

Note that the completent of $p_{a}(M_{{\rm free}} \cap \mu^{-1}(a)) = (M_{{\rm free}})_{a}$ has measure $0$ for the projection 
$p_{a} : \mu^{-1}(a) \longrightarrow M_{a}$ because the complement of $(M_{{\rm free}})_{a}$ is equal to the image of 
a finite union of submanifolds (or suborbifolds) of $\mu^{-1}(a)$ of codimension $\geq 2$. 
%However, the trivialization showed that the $T$-action locally does not depend on $\xi$ which shows that $p_{a}(M_{{\rm free}} \cap \mu^{-1}(a)) = M_{{\rm free}}_{a}$ 
%is equal to the complement of a finite union of submanifolds (or suborbifolds) of $M_{a}$ of codimension $\geq 2$. 
Since $p_{a}:M_{{\rm free}} \cap \mu^{-1}(a) \longrightarrow (M_{{\rm free}})_{a}$ is a principle $T$-fibration and ${\rm vol}(T) = 1$, we get that 
the volume of $M_{{\rm free}} \cap \mu^{-1}(a)$ is equal to the volume of $(M_{{\rm free}})_{a}$ with respect to the measure $dm_{a}$ 
induced by the reduced generalized Calabi-Yau structure on $M_{a}$. 
Because the complement of $(M_{{\rm free}})_{a}$ has measure 0, we have proved the formula. 
\end{proof}

\begin{rem}
For the density function $f$, one can show that $f$ is a piecewise polynomial of degree at most $n-l-k$. 
Moreover, in the case that $M$ is compact, the localization formula holds by applying the Atiyah-
Bott-Berline-Vergne localization theorem. Detailed statements and proofs can be seen in \cite{humi}. 
\end{rem}

%references

\end{document}